\documentstyle{article}
\begin{document}

\begin{center}
{\large \bf Hilbert H*-modules}\\
Michael Frank \copyright \\
frank@mathematik.uni-leipzig.d400.de or \\
frank@rz.uni-leipzig.de
\end{center}

\noindent
\rule{6in}{0.1mm}

\noindent
last update: 1.4.1996

\nocite{Biedenharn/Horwitz:84,Biedenharn/Horwitz:85}
\nocite{Delanghe/Brackx:78,Delanghe:79,Giellis:72,Goldstine/Horwitz:66}
\nocite{Molnar:92,Molnar:90,Molnar:92/2,Molnar:92/3,Molnar:93/1,Molnar:93/2}
\nocite{Molnar:91,Razon/Horwitz:92,Razon/Horwitz:91/1,Razon/Horwitz:91/2}
\nocite{Cabrera/Martinez/Rodriguez:95,Cabrera/Martinez/Rodriguez:89/2}
\nocite{Saworotnow:68,Saworotnow:70,Saworotnow:85,Sharma/Coulson:87}
\nocite{Soffer/Horwitz:83,Smith:74,Cnops:92}
\nocite{Saworotnow:80,Saworotnow:82,Smith:73,Bultheel:82,Kulkarni:94}
\nocite{Saworotnow:91,Cabrera:91,Saworotnow:81,Wiener/Masani:57}
\nocite{Wiener/Masani:58,Wiener/Masani:60,Saworotnow:83,Kokschal:95,Fuge:95}
\nocite{Masani:59,Masani:62,Masani:66,Truong-Van:81,Salehi:65,Salehi:66}
\nocite{Salehi:67,Rosenberg:64,Ellis/Gohberg/Lay:95,Ben-Artzi/Gohberg:94}
\nocite{Sen:82,Bultheel:80,Cabrera/Martinez/Rodriguez:92}
\nocite{Cabrera/ElMarrakchi/Martinez/Rodriguez:93,Razon/Horwitz/Biedenharn:89}

\newpage

\begin{center}
{\large \bf Hilbert modules over operator algebras}  \\
Michael Frank \copyright \\
frank@mathematik.uni-leipzig.d400.de or \\
frank@rz.uni-leipzig.de
\end{center}

\noindent
\rule{6in}{0.1mm}

\noindent
last update: 1.4.96

\vspace{1cm}

\medskip \noindent
The collection of sources below is far from being complete, but
nevertheless, instructive as guide into the research field.

{\small

\nocite{Blecher/Muhly/Paulsen:96,Blecher:96,Na:95,Muhly/Solel:95}
\nocite{Douglas/Paulsen:89,Blecher:95/1,Blecher:95/2,Blecher:95/3}
\nocite{Blecher/Muhly/Na:96/1,Blecher/Muhly/Paulsen:94,Muhly/Solel:93}
\nocite{Carlson/Clark:95,Misra/Pati:93,Douglas/Paulsen/Sah/Yan:95}
\nocite{Douglas:95,Carlson/Clark/Foias/Williams:94,Chen/Douglas:92}
\nocite{Carlson/Clark:93,Douglas:90,Sun/Li:88,Douglas:88,Douglas:86}

}

\end{document}